\documentstyle[graphicx]{article}[15pt]
\input amssym.def

\catcode`@=12
 \long\def\@makefntext#1{\noindent #1}
\newskip\tabcentering \tabcentering=1000pt plus 1000pt minus 1000pt

\def\MCH#1#2{\setbox0=\hbox{\raise#1\hbox{#2}}\smash{\box0}}

\def\@evenfoot{}\def\@oddfoot{}

\def\sec#1{\vspace{5mm}\leftline{\bf #1}\vspace{3mm}}

\catcode`@=12

\def\bc{\begin{center}}
\def\ec{\end{center}}

\def\hang{\hangindent\parindent}
\def\textindent#1{\indent\llap{\qquad #1\ \ \enspace}\ignorespaces}
\def\ref{\par\hang\textindent}

\def\a1{(a_1, a_2, \cdots, a_n)}

\def\a{\alpha}

\begin{document}
\thispagestyle{empty}
\vspace*{-3.0truecm}
\noindent
\vspace{1 true cm}
 \bc{\large\bf Concentration of mass in the   pressureless limit of   the  Euler equations of one-dimensional compressible
fluid flow

\footnotetext{$^{*}$Corresponding author. \\
\indent \,\,\,\,\,\,\,\,E-mail address:  zqshao@fzu.edu.cn.\\
 }}\ec

 \vspace*{0.2 true cm}
\bc{\bf Shouqiong Sheng$^{a}$, Zhiqiang  Shao$^{a, *}$   \\
{\it $^{a}$College of Mathematics and Computer Science, Fuzhou University, Fuzhou 350108, China}
 }\ec

 \vspace*{2.5 true mm}
\setlength{\unitlength}{1cm}
\begin{picture}(20,0.1)
\put(-0.6,0){\line(1,0){14.5}}
\end{picture}

 \vspace*{2.5 true mm}
\noindent{\small {\small\bf Abstract}

 \vspace*{2.5 true mm}In this paper, we study the limiting behavior of Riemann solutions to the  Euler equations of one-dimensional
 compressible fluid flow as $\gamma$ tends to one. We show that the limit solution forms the delta wave to the pressureless Euler system  of one-dimensional
 compressible fluid flow in the distribution sense.
Some numerical results
exhibiting the phenomena of concentration are also presented.

 \vspace*{2.5 true mm}
\noindent{\small {\small\bf MSC: } 35L65;  35L67

 \vspace*{2.5 true mm}
\noindent{\small {\small\bf Keywords:} Riemann solutions;  Delta  wave;  Pressureless limit; Euler equations of one-dimensional compressible
fluid flow; Numerical simulations

 \vspace*{2.5 true mm}
\setlength{\unitlength}{1cm}
\begin{picture}(20,0.1)
\put(-0.6,0){\line(1,0){14.5}}
\end{picture}



\baselineskip 15pt
 \sec{\Large\bf 1.\quad  Introduction }
  The Euler equations of one-dimensional compressible fluid flow  read (cf. [8]):
$$ \left\{\begin{array}{ll} \rho_{t}+(\rho u)_{x}=0,\\u_{t}+(\frac{u^{2}}{2}+p (\rho))_{x}=0,\end{array}\right .\eqno{(1.1)}
$$
where the nonlinear   function  $p(\rho)=\frac{\theta}{2}\rho^{\gamma-1},$   $\theta=\frac{\gamma-1}{2}$  and  $\gamma\in(1,2)$   is a constant.

System (1.1) was firstly derived by Earnshaw [8] in 1858 for isentropic flow and is also viewed as
 the Euler equations of one-dimensional compressible fluid flow [11]. where $\rho$ denotes the density, $u$ the velocity, and $p(\rho)$ the
pressure of the fluid.
System (1.1) has other different physical backgrounds. For instance, it is a scaling
limit system of Newtonian dynamics with long-range interaction for a continuous
distribution of mass   in $R$ [15, 16] and also a hydrodynamic
limit for the Vlasov equation [1].

The  solutions for system (1.1)  were widely studied by many scholars (see [5-7,  13-14, 17, 23] ). In particular,  the existence of global weak solutions of the Cauchy problem  was first established by DiPerna [7] for the case of $1 <\gamma<3$ by using the Glimm's
scheme method.
 Using the result of DiPerna [7], Li [13] obtained a global weak solution to the Cauchy problem for
the case $  -1 < \gamma < 1.$  Using the theory of compensated compactness coupled with some
basic ideas of the kinetic formulation, Lu [14]  established an existence theorem for global entropy solutions for the case $\gamma > 3$. Cheng [6] also used the same methods as in [14]
 to obtain  the existence of global entropy solutions for the Cauchy
problem with a uniform amplitude bound for
the case $1 <\gamma  < 3.$

In this paper, we are interested in
   the Riemann problem for (1.1)   with initial data
$$ (\rho, u)(0, x) =\left\{\begin{array}{ll} (\rho_{-},
u_{-}),\,\,\,\,x< 0,\\(\rho_{+},
u_{+}),\,\,\,\,x> 0,\end{array} \right.\eqno{(1.2)}$$
where  $\rho_{\pm}>0$ and $u_{\pm}$   are given constant states. We assume that
$u_{+}<u_{-}$,  and $\gamma\in(1,2)$.

  System (1.1) is just like a hyperbolic system for conservation laws of the form
  $$ \partial_{t}U+\partial_{x}F(U)=0,\eqno{(1.3)} $$
with
$$ U=\left(\begin{array}{c}\rho \\u
 \end{array}\right), ~~F(U)=\left(\begin{array}{c}\rho u \\\frac{u^{2}}{2}+\frac{\gamma-1}{4}\rho^{\gamma-1}
 \end{array}\right)=0.$$

  When $\gamma\rightarrow 1$, the limiting system of (1.1) formally becomes  the pressureless Euler system of one-dimensional compressible fluid flow,
  $$ \left\{\begin{array}{ll} \rho_{t}+(\rho u)_{x}=0,\\u_{t}+(\frac{u^{2}}{2})_{x}=0.\end{array}\right .\eqno{(1.4)}
$$

 Let us turn to the Euler system  of power law in Eulerian coordinates,
$$ \left\{\begin{array}{ll} \rho_{t}+(\rho u)_{x}=0,\\(\rho u)_{t}+
 (\rho u^{2}+p(\rho) )_{x}=0,\end{array}\right .
\eqno{(1.5)}$$  When the pressure tends to zero or a constant, the Euler system  (1.5)  formally tends to the zero pressure gas dynamics.  In earlier seminal papers,
Chen and Liu   [3] first showed the formation of $\delta$-shocks and vacuum states of the Riemann solutions to the Euler system (1.5) for polytropic gas by taking limit $\varepsilon \rightarrow 0+$ in the model $p(\rho)=\varepsilon\rho^{\gamma}/\gamma ( \gamma >1)$,
which describe the phenomenon of concentration and cavitation rigorously in mathematics.  Further, they also obtained the same results for the Euler
equations for nonisentropic
fuids in [4]. The same problem for the  Euler equations  (1.5) for
isothermal case $( \gamma =1)$ was studied by Li [12], in which he proved that when temperature drops to zero, the
solution containing two shock waves converges to the delta shock solution to the transport equations
and the solution containing two rarefaction waves converges to the solution involving vacuum to the
transport equations.  Recently,  Muhammad Ibrahim, Fujun Liu and Song Liu
[9] showed the same phenomenon of concentration also exists in the mode $p(\rho) = \rho^{\gamma}$ $(0 < \gamma < 1)$ as $\gamma \rightarrow 0$, which is the case that the pressure goes to a constant. Namely, they  showed rigorously the formation of delta wave with the limiting behavior of Riemann solutions to the Euler equations (1.5).

    Motivated by [9], for the  Euler system (1.1) of one-dimensional compressible fluid flow, we show the same phenomenon of concentration also exists in the case $1 < \gamma < 2$ as $\gamma \rightarrow1$. To our knowledge,   the Riemann problem for  system (1.1) when $1 < \gamma < 2$ was not studied before in the
literature.

 The Riemann problem for  system (1.1) when  $\gamma=-1$ was solved by Cheng et al. [5], where the delta shock wave solutions were obtained. For the case  $\gamma=-1$, system (1.1)  can also be called the Chaplygin Euler equations of compressible fluid flow. The Chaplygin gas was
introduced   by  Chaplygin [2], Tsien [21] and von Karman [22] as  a suitable mathematical approximation for calculating the lifting force on a wing of an airplane in aerodynamics.
 Zhang et al. [23] also showed Interaction of delta shock waves for the Chaplygin Euler
equations of compressible fluid flow with split delta
functions in this case, and their numerical simulations completely coinciding with the theoretical analysis were also exhibited. Recently, the Riemann problem for the Euler equations of  compressible fluid flow with the
generalized Chaplygin gas was studied by Pang et al. [17].
\vskip 0.1in
The paper is organized as follows.
 In Section 2, we give some preliminaries on the delta wave to the pressureless Euler system
  of one-dimensional compressible fluid flow. In Section 3, we display some results on  the Riemann solutions of (1.1)  with initial data (1.2) when $1 < \gamma < 2$.
In Section 4, we show rigorously the formation of delta wave with the limiting behavior of Riemann solutions to the Euler system
(1.1).
In section 5, we
present some representative numerical results to examine the phenomenon of
concentration and the formation process of the delta wave in  Riemann solutions to the Euler system  of one-dimensional compressible fluid flow
 as the adiabatic exponent  $\gamma$  decreases.

\baselineskip 15pt
 \sec{\Large\bf 2.\quad   Preliminaries }
  In this section, let us  briefly review the delta wave  to the Riemann problem
for the  pressureless Euler system (1.4). As mentioned in [9],
 the Riemann solutions to the pressureless Euler system (1.4) contain delta wave, in which the solution is a delta function supposed on a smooth curve,
see, such as, Sheng and  Zhang [19], Joseph [10], Shen et al. [18], Sun [20].

We first introduce  the definition of the two-dimensional weighted delta function  as  follows.

 \vskip 0.1in
  \noindent{\small {\small\bf Definition 2.1.}
 A two-dimensional weighted delta function $w(s)\delta_{S}$ supported on a smooth curve $ S=\{(t(s),x(s)):a<s<b\}$ is defined by
$$ \langle w(t)\delta_{S},\varphi(t,x)\rangle=\int_{a}^{b}w(t(s))\varphi(t(s),x(s))ds, \eqno{(2.1)}$$
for all test functions $\varphi(t,x) \in C_{0}^{\infty}([0,+\infty)\times (-\infty, +\infty)).$ \\

For the Riemann problem  with $u_{+}<u_{-}$, we can construct a dirac-measured solution with
 parameter $\sigma$  as follows,
 $$  \rho(t,x)=\rho_{0}(t,x)+w(t)\delta_{S}, ~~u(t,x)=u_{0}(t,x), \eqno{(2.2)}
 $$
 where $ S=\{(t,\sigma t):0\leq t<+\infty\} $,
 $$\rho_{0}(t,x)=\left\{
                  \begin{array}{ll}
                    \rho_{-}, & \hbox{$x<\sigma t $,} \\
                    \rho_{+}, & \hbox{$x>\sigma t$,}
                  \end{array}
                \right.
\eqno{}$$$$u_{0}(t,x)=\left\{
                  \begin{array}{ll}
                   u_{-}, & \hbox{$x<\sigma t $,} \\
                    \sigma, & \hbox{$x=\sigma t$,} \\
                    u_{+}, & \hbox{$x>\sigma t$,}
                  \end{array}
                \right.
\eqno{}$$
and$$ w(t)=t(\sigma[\rho]-[\rho u]),\eqno{(2.3)}$$
where $[q]=q_{+} -q_{-}$ denotes the jump of  function $q$  across the discontinuity. The dirac-measured
solution $(\rho, u)$ constructed above is  known as ``delta wave" to the  pressureless Euler system (1.4) if
$$ \langle \rho,\varphi_{t}\rangle+\langle \rho u,\varphi_{x}\rangle=0,\eqno {(2.4)}
$$$$ \langle u,\varphi_{t}\rangle+\langle \frac{u^{2}}{2}, \varphi_{x}\rangle=0,\eqno {(2.5)}
$$
hold  for any test function $\varphi (t,x)\in C_{0}^{\infty}([0,+\infty)\times (-\infty, +\infty))$, where
$$ \langle \rho, \varphi\rangle=
       \int_{0}^{+\infty}\int_{-\infty}^{+\infty}\rho_{0}(t,x)\varphi (t,x)dxdt+\langle w(t)\delta_{S}, \varphi(t,x)\rangle,$$
       $$\langle \rho u, \varphi\rangle=
      \int_{0}^{+\infty}\int_{-\infty}^{+\infty}\rho_{0}(t,x)u_{0}(t,x)\varphi (t,x)dxdt+\langle \sigma  w(t)\delta_{S}, \varphi(t,x)\rangle,
  $$
and $u$ has the similar integral identity as above.  Then the following generalized Rankine-Hugoniot conditions
$$
\left\{
     \begin{array}{ll}
       \frac{dx}{dt}=\sigma, \\
       \frac{dw(t)}{dt}=\sigma [\rho]-[\rho u], \\
       \sigma [u]=[\frac{1}{2} u^{2}],
     \end{array}
   \right.\eqno (2.6)
$$
hold, where $[\rho]= \rho_{+}-\rho_{-}$,   with initial data
$$(x, w)(0) = (0,  0).\eqno{(2.7)}$$

Solving the system of the simple ordinary differential equations (2.6) with initial data (2.7),
we have$$
 w(t)=\frac{1}{2}(\rho_{-}+\rho_{+})(u_{-}-u_{+})\,t, ~~~~\sigma=\frac{1}{2}(u_{-}+u_{+})$$
 fulfilling the entropy condition $u_{+}<\sigma<u_{-}$, see [18].
\vskip 0.1in
\noindent{\small {\small\bf Remark 2.1.} The entropy condition $u_{+}<\sigma<u_{-}$ means that all characteristic lines  on both sides of the discontinuity are incoming. So it is a
overcompressive condition.

\baselineskip 15pt
 \sec{\Large\bf 3.\quad   Riemann problem   for   Euler
equations
 }
{\Large\bf \quad of one-dimensional compressible fluid flow
}

 \indent

 In this section, we present some preliminary knowledge for system (1.1) and construct the
Riemann solutions of (1.1) with initial data (1.2).

 The eigenvalues of  system (1.1) are
 $$ \lambda_{1}=u-\frac{\gamma-1}{2}\rho^{\frac{\gamma-1}{2}},\,\,\,\,\,\,\,\lambda_{2}=u+\frac{\gamma-1}{2}\rho^{\frac{\gamma-1}{2}},\eqno{(3.1)}
$$
with the corresponding right eigenvectors
$$\overrightarrow{r}_{1} =(1, -\frac{\gamma-1}{2}\rho^{\frac{\gamma-3}{2}})^{T}, \,\,
\overrightarrow{r}_{2} =(1, \frac{\gamma-1}{2}\rho^{\frac{\gamma-3}{2}})^{T},
 \eqno{(3.2)}
$$satisfying
$$\nabla\lambda_{1}\cdot \overrightarrow{r}_{1}=-\frac{(\gamma-1)(\gamma+1)}{4}\rho^{\frac{\gamma-3}{2}}<0,$$ $$\nabla\lambda_{2}\cdot \overrightarrow{r}_{2}=\frac{(\gamma-1)(\gamma+1)}{4}\rho^{\frac{\gamma-3}{2}}>0.$$Therefore,  system (1.1)  is strictly hyperbolic for $\rho>0$,
 both   characteristic fields  are genuinely nonlinear  and the associated waves are either shock waves  or rarefaction waves.

Since system (1.1)  and the Riemann data (1.2) are invariant under stretching of coordinates:
$(t, x)\rightarrow (\kappa t, \kappa x)~(\kappa$ is a constant),  we seek the self-similar solution $$(\rho,u)(t, x)=( \rho,u)(\xi),\,\,\,\,\xi=\frac{x}{t}.$$
 Then  Riemann problem (1.1) and (1.2) is reduced to the following boundary value problem of  ordinary differential equations:
$$ \left\{\begin{array}{ll}
   -\xi\rho_{\xi}+(\rho u)_{\xi}=0,\\
   -\xi u_{\xi}+\big(\frac{u^{2}}{2}+\frac{\gamma-1}{4}\rho^{\gamma-1} \big)_{\xi}=0,\end{array}\right .\eqno{(3.3)}
$$
with $(\rho,u)(\pm\infty)=( \rho_{\pm},u_{\pm}).$\\
\indent
For any smooth solution, system (3.3) can be rewritten as

$$\left(
  \begin{array}{ccc}u-\xi&\rho
    \\
    \frac{(\gamma-1)^{2}}{4}+\rho^{\gamma-2}&u-\xi
  \end{array}\right)\left(\begin{array}{cccc} \rho_{\xi}\\u_{\xi}
 \end{array}\right)=0.\eqno{(3.4)}$$
It provides either the general solution (constant state)
$$( \rho, u)(\xi)={\mathrm constant }   \,\,(\rho>0), $$
or  the
1-rarefaction wave
$$R_{1}(\rho_{-},u_{-}):\,\,\left\{\begin{array}{ll} \xi=\lambda_{1}=u-\frac{\gamma-1}{2}\rho^{\frac{\gamma-1}{2}},
\\u-u_{-}=-(\rho^{\frac{\gamma-1}{2}}-\rho_{-}^{\frac{\gamma-1}{2}}),
\,\,\,\,\rho<\rho_{-}, u>u_{-},\end{array} \right.\eqno{(3.5)}$$  or  the 2-rarefaction wave
$$R_{2}(\rho_{-},u_{-}):\,\,\left\{\begin{array}{ll} \xi=\lambda_{2}=u+\frac{\gamma-1}{2}\rho^{\frac{\gamma-1}{2}},
\\u-u_{-}=\rho^{\frac{\gamma-1}{2}}-\rho_{-}^{\frac{\gamma-1}{2}},
\,\,\,\,\rho>\rho_{-}, u>u_{-}.\end{array} \right.\eqno{(3.6)}$$

Differentiating
 the second equation of  (3.5) with respect to $\rho$ yields $$u_{\rho} = -\frac{\gamma-1}{2}\rho^{\frac{\gamma-3}{2}}<0,$$
 and
subsequently,
$$u_{\rho\rho} =-\frac{(\gamma-1)(\gamma-3)}{4}\rho^{\frac{\gamma-5}{2}}>0,$$
which mean that the 1-rarefaction wave curve $R_{1}(\rho_{-},u_{-})$ is monotonic decreasing and convex in the $(\rho, u)$ phase plane
$(\rho> 0)$. Similarly, from the second equation of (3.6), we have $u_{\rho}> 0$ and $u_{\rho\rho} < 0,$ which mean that
the 2-rarefaction wave curve $R_{2}(\rho_{-},u_{-})$ is monotonic increasing and concave in the $(\rho, u)$ phase plane
$(\rho> 0)$.  Moreover, it  can be concluded  from (3.5) that
 $\lim\limits_{\rho\rightarrow 0^{+}}u= u_{-}+\rho_{-}^{\frac{\gamma-1}{2}}$  for the 1-rarefaction wave curve $R_{1}(\rho_{-},u_{-})$, which implies that
$R_{1}(\rho_{-},u_{-})$ intersects the $u$-axis at the point $(0,  \widetilde{u}_{\ast})$,  where $\widetilde{u}_{\ast}$ is determined by
$\widetilde{u}_{\ast}= u_{-}+\rho_{-}^{\frac{\gamma-1}{2}}$. It can also be  seen from (3.6) that $\lim\limits_{\rho\rightarrow +\infty}u=+\infty$
 for the 2-rarefaction
wave curve $R_{2}(\rho_{-},u_{-})$.

  Let $x=\sigma t$ be a discontinuity  of bounded discontinuous  solutions of (1.1), the following Rankine-Hugoniot conditions
$$\left\{
    \begin{array}{ll}
      \sigma[\rho]=[\rho u], \\
     \sigma[ u]=[\frac{u^{2}}{2}+\frac{\gamma-1}{4}\rho^{\gamma-1}],
    \end{array}
  \right.\eqno (3.7)
$$
hold, where $[\rho]=\rho-\rho_{-}$, etc. From (3.7) we can get
$$
      \sigma=\frac{[\rho u]}{[\rho]},$$
     $$ u=u_{-}\pm(\rho-\rho_{-})\sqrt{\frac{\frac{\gamma-1}{2}[\rho^{\gamma-1}]}{(\rho+\rho_{-})[\rho]}},
     \eqno (3.8)$$
where $\sigma$, $(\rho_{-},u_{-})$ and $(\rho_{},u_{})$ are the shock speed, the left state and the right state, respectively.

1-$shock$  $curve$  $S_{1}(\rho_{-}, u_{-})$:

The classical Lax
 entropy conditions   imply  that the propagation speed $\sigma$ for the 1-shock wave
has to be satisfied with$$\sigma<\lambda_{1}(\rho_{-},u_{-}),\,\,\,\,\lambda_{1}(\rho,u)<\sigma<\lambda_{2}(\rho,u).\eqno (3.9)
$$
On the other hand, from the first equation of (3.7), we have
$$\sigma=\frac{\rho u -\rho_{-}u_{-}}{\rho-\rho_{-}} = u_{-}+\frac{\rho}{\rho-\rho_{-}}(u-u_{-}).$$
Thus,  it follows from the  first inequality  of (3.9) that
$$\frac{\rho}{\rho-\rho_{-}}(u-u_{-})<-\frac{\gamma-1}{2}\rho_{-}^{\frac{\gamma-1}{2}}<0,$$
which means that $u - u_{-}$ and $\rho- \rho_{-}$ have different signs. Then
from (3.8)  we have $$u=u_{-}-\sqrt{\frac{\frac{\gamma-1}{2}(\rho^{\gamma-1}-\rho_{-}^{\gamma-1})}{(\rho+\rho_{-})(\rho-\rho_{-})}}(\rho-\rho_{-}).
$$
 If $u>u_{-}$, then $\rho<\rho_{-}$, and
 $$\sigma- u_{-}=\frac{\rho}{\rho-\rho_{-}}(u-u_{-})=-\rho\sqrt{\frac{\frac{\gamma-1}{2}(\rho^{\gamma-1}-\rho_{-}^{\gamma-1})}{(\rho+\rho_{-})(\rho-\rho_{-})}}=-\frac{\gamma-1}{2}
 \overline{\rho}^{\frac{\gamma-2}{2}}\rho\sqrt{\frac{2}{\rho+\rho_{-}}},
$$
for some $\bar{\rho}\in(\rho,\rho_{-}).$  By direct calculation, we have
$$\frac{\gamma-1}{2}\rho_{-}^{\frac{\gamma-1}{2}}-\frac{\gamma-1}{2}
 \overline{\rho}^{\frac{\gamma-2}{2}}\rho\sqrt{\frac{2}{\rho+\rho_{-}}}>\frac{\gamma-1}{2}\bigg(\rho_{-}^{\frac{\gamma-1}{2}}-
 \rho^{\frac{\gamma-2}{2}}\rho\sqrt{\frac{2}{\rho+\rho_{-}}}\bigg)>\frac{\gamma-1}{2}(\rho_{-}^{\frac{\gamma-1}{2}}-\rho^{\frac{\gamma-1}{2}})>0,
$$ which implies that
$$\sigma-u_{-} >-\frac{\gamma-1}{2}\rho_{-}^{\frac{\gamma-1}{2}}.
$$
This  contradicts with $\sigma<\lambda_{1}(\rho_{-},u_{-})$.  Hence we get the 1-shock wave curve $S_{1}(\rho_{-},u_{-})$
  in the phase plane,
$$u=u_{-}-\sqrt{\frac{\frac{\gamma-1}{2}(\rho^{\gamma-1}-\rho_{-}^{\gamma-1})}{(\rho+\rho_{-})(\rho-\rho_{-})}}(\rho-\rho_{-}),~~~\rho>\rho_{-},u<u_{-}.\eqno{(3.10)}
$$

2-$shock~curve~S_{2}(\rho_{-},u_{-})$:

Similarly, the propagation speed $\sigma$ for the 2-shock wave should satisfy
$$\lambda_{1}(\rho_{-},u_{-})<\sigma<\lambda_{2}(\rho_{-},u_{-}),\,\,\,\,\lambda_{2}(\rho,u)<\sigma.
$$Then,  we can get the 2-shock wave curve $S_{2}(\rho_{-},u_{-})$ in the phase plane,
 $$u=u_{-}+\sqrt{\frac{\frac{\gamma-1}{2}(\rho^{\gamma-1}-\rho_{-}^{\gamma-1})}{(\rho+\rho_{-})(\rho-\rho_{-})}}(\rho-\rho_{-}),~~  \rho<\rho_{-},u<u_{-}.\eqno{(3.11)}
 $$

Differentiating
 both sides of (3.10) with respect to $\rho$ gives that for $\rho>\rho_{-}$,
$$u_{\rho}=-\frac{1}{2}\sqrt{\frac{\frac{\gamma-1}{2}(\rho+\rho_{-})}{(\rho^{\gamma-1}-\rho_{-}^{\gamma-1})(\rho-\rho_{-})}}
\frac{(\gamma-1)\rho^{\gamma-2}(\rho-\rho_{-})(\rho+\rho_{-})+2\rho_{-}(\rho^{\gamma-1}-\rho_{-}^{\gamma-1})}{(\rho+\rho_{-})^{2}}<0,$$
which means that the 1-shock wave curve $S_{1}(\rho_{-},u_{-})$ is monotonic decreasing in the $(\rho, u) $ phase plane  ($\rho>\rho_{-})$. Similarly,
from  (3.11), for $\rho<\rho_{-}$ we have $u_{\rho} > 0, $  which means that the 2-shock wave
curve $S_{2}(\rho_{-},u_{-})$ is monotonic increasing in the $(\rho, u) $ phase plane  ($\rho<\rho_{-})$. In addition, It can be seen from (3.11)
that $\lim\limits_{\rho\rightarrow 0^{+}}u= u_{-}-\sqrt{\frac{\gamma-1}{2}}\rho_{-}^{\frac{\gamma-1}{2}}$ for the 2-shock wave curve $S_{2}(\rho_{-},u_{-})$, which implies that $S_{2}(\rho_{-},u_{-})$ intersects the $u$-axis at the point $(0,  \widetilde{u}_{\ast\ast} )$,  where $\widetilde{u}_{\ast\ast}$
is determined by  $\widetilde{u}_{\ast\ast}=u_{-}-\sqrt{\frac{\gamma-1}{2}}\rho_{-}^{\frac{\gamma-1}{2}}.$
 It can also be derived from (3.10) that $\lim\limits_{\rho\rightarrow +\infty}u= -\infty$
 for the 1-shock wave curve $S_{1}(\rho_{-},u_{-})$.

In the $(\rho, u) $ phase plane, through a given point $(\rho_{-}, u_{-})$, we draw the elementary wave curves $R_{j}$ and $S_{j}$ (j=1,2).  We find that the elementary wave curves divide the  $(\rho, u)$ phase plane into five regions (see Fig. 1).  According to the right
state $(\rho_{+},u_{+})$ in the different regions, one can construct the unique global solution to the Riemann problem (1.1) and (1.2) as follows:

(1) $(\rho_{+},u_{+})\in I(\rho_{-},u_{-}):$ $R_{1}+R_{2};$

 (2)$(\rho_{+},u_{+})\in II(\rho_{-},u_{-}):$ $S_{1}+R_{2};$

 (3)$(\rho_{+},u_{+})\in III(\rho_{-},u_{-}):$ $R_{1}+S_{2};$

 (4)$(\rho_{+},u_{+})\in IV(\rho_{-},u_{-}):$ $S_{1}+S_{2};$

  (5)$(\rho_{+},u_{+})\in V(\rho_{-},u_{-}):$ $R_{1}+\mathrm{Vac}+R_{2},$\\
 where   ``+" means ``followed by".

\hspace{65mm}\setlength{\unitlength}{0.8mm}\begin{picture}(80,66)
\put(-50,18){\vector(0,2){35}}
 \put(-48,0){\vector(2,0){103}}  \put(-53,49){$\rho$}
\put(56,-1){$u$}
\put(-38,7){$S_{2}$}\put(6,5){}
\put(-32,48){$S_{1}$}\put(32,48){$R_{2}$}\put(42,19){$R_{2}$}\put(12,5){$R_{1}$}\put(45,6){V }\put(25,-5){$\widetilde{u}_{\ast}$ }\put(-40,-5){$\widetilde{u}_{\ast\ast}$ }
\put(1,12){$(\rho_{-},
u_{-})$}
\put(-5,6){III }\put(44,6){}\put(-5,29){II}
\put(-36,20){IV}\put(26,20){I}
\qbezier(-40,0)(14,12)(42,52)\qbezier(25,0)(-14,12)(-42,52)\qbezier(25,0)(36,06)(48,18)
\end{picture}
\vspace{0.6mm}  \vskip 0.2in \centerline{\bf Fig. 1.\,\,    Curves of  elementary waves.
   } \vskip 0.1in \indent

 We are interested here in  the case $S_{1}+S_{2}$ that there exists a unique intermediate state $(\rho_{\ast},u_{\ast})$
 such that $(\rho_{\ast},u_{\ast})\in S_{1}(\rho_{-},u_{-})$ and $(\rho_{+},u_{+})\in S_{2}(\rho_{\ast},u_{\ast})$, i.e.,
$$ u_{\ast}=u_{-}-\sqrt{\frac{\frac{\gamma-1}{2}(\rho_{\ast}^{\gamma-1}-\rho_{-}^{\gamma-1})}{(\rho_{\ast}+\rho_{-})(\rho_{\ast}-\rho_{-})}}(\rho_{\ast}-\rho_{-}),
~~\rho_{\ast}>\rho_{-},~u_{\ast}<u_{-},\eqno{(3.12)}  $$
$$ u_{+}=u_{\ast}+\sqrt{\frac{\frac{\gamma-1}{2}(\rho_{+}^{\gamma-1}-\rho_{\ast}^{\gamma-1})}{(\rho_{\ast}+\rho_{+})(\rho_{+}-\rho_{\ast})}}(\rho_{+}-\rho_{\ast}),
    ~~\rho_{+}<\rho_{\ast},~u_{+}<u_{\ast}, \eqno{(3.13)}  $$
with the shock speed:
$$ \sigma_{1}=\frac{\rho_{\ast}u_{\ast}-\rho_{-}u_{-}}{\rho_{\ast}-\rho_{-}},
 ~~\sigma_{2}=\frac{\rho_{+}u_{+}-\rho_{\ast}u_{\ast}}{\rho_{+}-\rho_{\ast}},\eqno{(3.14)}
$$
respectively. In this case, the Riemann solution is
$$ (\rho, u)(t, x) =\left\{\begin{array}{ll} (\rho_{-},
u_{-}),\,\,\,\,x< \sigma_{1}t,\\(\rho_{\ast},
u_{\ast}),\,\,\,\,\sigma_{1}t< x < \sigma_{2}t,\\(\rho_{+},
u_{+}),\,\,\,\,x> \sigma_{2}t.\end{array} \right.\eqno{(3.15)}$$

\baselineskip 15pt
 \sec{\Large\bf 4.\quad   Formation of delta wave}In this section, we  study the limiting behavior of   solutions of the Riemann problem (1.1)  and (1.2) with the assumption $u_{+}<u_{-}$ as $\gamma$ goes to one.
 Then we show the limit is the delta wave  of  the pressureless Euler system (1.4).

 \vskip 0.1in
\noindent{\small {\small\bf Lemma 4.1.} If $u_{+}<u_{-}$, then there is a sufficiently small $\gamma_{0} > 0$ such that
$(\rho_{+}, u_{+})\in IV( \rho_{-}, u_{-})$ as $1<\gamma <1 +\gamma_{0}$.

\vskip 0.1in
\noindent{\small {\small\bf Proof.} If $\rho_{+} = \rho_{-}$, then $(\rho_{+}, u_{+})\in IV( \rho_{-}, u_{-})$  for any $\gamma\in (1,2)$. Thus, we only need to consider the
case $\rho_{+} \neq\rho_{-}$.

It can be derived from (3.10) and (3.11) that all possible states $(\rho, u)$ that can be connected to the
left state $(\rho_{-}, u_{-})$ on the right by a 1-shock wave $S_{1}$ or a 2-shock wave $S_{2}$ should satisfy
$$   S_{1}: \,\,\,\,u=u_{-}-\sqrt{\frac{\frac{\gamma-1}{2}(\rho^{\gamma-1}-\rho_{-}^{\gamma-1})}{(\rho+\rho_{-})(\rho-\rho_{-})}}(\rho-\rho_{-}),~~~\rho>\rho_{-},\eqno{(4.1)}
$$
 $$S_{2}: \,\,\,\,u=u_{-}+\sqrt{\frac{\frac{\gamma-1}{2}(\rho^{\gamma-1}-\rho_{-}^{\gamma-1})}{(\rho+\rho_{-})(\rho-\rho_{-})}}(\rho-\rho_{-}),~~  \rho<\rho_{-}.\eqno{(4.2)}
 $$
If $\rho_{+} \neq\rho_{-}$ and $(\rho_{+}, u_{+})\in IV( \rho_{-}, u_{-})$, then from Fig. 1, (4.1) and (4.2), we have$$  u_{+}<u_{-}-\sqrt{\frac{\frac{\gamma-1}{2}(\rho_{+}^{\gamma-1}-\rho_{-}^{\gamma-1})}{(\rho_{+}+\rho_{-})(\rho_{+}-\rho_{-})}}(\rho_{+}-\rho_{-}),~~~\mathrm {}  \,\,\rho_{+}>\rho_{-},\eqno{(4.3)}
$$
$$  u_{+}<u_{-}+\sqrt{\frac{\frac{\gamma-1}{2}(\rho_{+}^{\gamma-1}-\rho_{-}^{\gamma-1})}{(\rho_{+}+\rho_{-})(\rho_{+}-\rho_{-})}}(\rho_{+}-\rho_{-}),~~~\mathrm {}  \,\,\rho_{+}<\rho_{-},\eqno{(4.4)}
$$
which implies that
$$ \sqrt{\frac{\frac{\gamma-1}{2}(\rho_{+}^{\gamma-1}-\rho_{-}^{\gamma-1})}{\rho_{+}^{2}-\rho_{-}^{2}}}< \frac{u_{-}-u_{+}}{|\rho_{+}-\rho_{-}|}. \eqno{(4.5)}$$
Since
$$  \lim_{{\gamma\rightarrow1}}\sqrt{\frac{\frac{\gamma-1}{2}(\rho_{+}^{\gamma-1}-\rho_{-}^{\gamma-1})}{\rho_{+}^{2}-\rho_{-}^{2}}}=0,
  \eqno{(4.6)}$$
it follows that there exists $\gamma_{0}>0$ small enough such that, when $ 1<\gamma <1 +\gamma_{0}$,  we have
$$ \sqrt{\frac{\frac{\gamma-1}{2}(\rho_{+}^{\gamma-1}-\rho_{-}^{\gamma-1})}{\rho_{+}^{2}-\rho_{-}^{2}}}< \frac{u_{-}-u_{+}}{|\rho_{+}-\rho_{-}|}. \eqno{}$$
  Then,  it is obvious
that $(\rho_{+}, u_{+})\in IV( \rho_{-}, u_{-})$ when $ 1<\gamma <1 +\gamma_{0}$.
 The proof is completed. $~\Box$

\vskip 0.1in

When  $1<\gamma <1+\gamma_{0}$, namely $(\rho_{+}, u_{+})\in IV( \rho_{-}, u_{-})$, suppose that $(\rho_{\ast},
u_{\ast})$ is the intermediate state
connected with $(\rho_{-}, u_{-})$ by a 1-shock wave $S_{1}$ with the speed $\sigma_{1}$, and $(\rho_{+}, u_{+})$ by a 2-shock wave $S_{2}$  with the speed $\sigma_{2},$
then it follows from
 (3.9) and (3.10) that
$$ u_{-}-u_{+}=\sqrt{\frac{\frac{\gamma-1}{2}(\rho_{\ast}^{\gamma-1}-\rho_{-}^{\gamma-1})}{(\rho_{\ast}+\rho_{-})(\rho_{\ast}-\rho_{-})}}(\rho_{\ast}-\rho_{-})+
\sqrt{\frac{\frac{\gamma-1}{2}(\rho_{+}^{\gamma-1}-\rho_{\ast}^{\gamma-1})}{(\rho_{\ast}+\rho_{+})(\rho_{+}-\rho_{\ast})}}(\rho_{\ast}-\rho_{+}), \,\,\,\,   \rho_{\ast}>\rho_{\pm}.\eqno{(4.7)}
$$

 Then we have the following lemma.

 \vskip 0.1in
\noindent{\small {\small\bf Lemma 4.2.}
$\lim\limits_{\gamma\rightarrow1}\rho_{\ast}=+\infty,$ and $ \lim\limits_{\gamma\rightarrow1}\frac{\gamma-1}{2}\rho_{\ast}^{\gamma-1}=:a
=\frac{(u_{-}-u_{+})^{2}}{4}$.

\vskip 0.1in
\noindent{\small {\small\bf Proof.} Let $ \lim\limits_{\gamma\rightarrow1}\inf\rho_{\ast}=\alpha$, and $\lim\limits_{\gamma\rightarrow1}\sup\rho_{\ast}=\beta$.

If $ \alpha<\beta$ , then by the continuity of $\rho_{\ast}(\gamma)$, there exists a sequence $  \{\gamma_{n}\}_{n=1}^{\infty}\subseteq(1,2)$
such that
$$ \lim_{n\rightarrow +\infty}\gamma_{n}=1,\,\,\mathrm{ and}\,\, \lim_{n\rightarrow +\infty}\rho_{\ast}(\gamma_{n})=c,$$
for  some $ c\in(\alpha,\beta).$ Then substituting the sequence  into the right hand side of (4.7),  and  taking the limit $n\rightarrow +\infty$,
 we have
$$ \lim_{n\rightarrow+\infty}\sqrt{\frac{\frac{\gamma_{n}-1}{2}(\rho_{\ast}(\gamma_{n})^{\gamma_{n}-1}-\rho_{\pm}^{\gamma_{n}-1})}
{\rho_{\ast}^{2}(\gamma_{n})-\rho^{2}_{\pm}}}(\rho_{\ast}(\gamma_{n})-\rho_{\pm})=0.
\eqno{(4.8)} $$
Thus,  we can obtain from (4.7) that
$$ u_{-}-u_{+}=0,$$
which contradicts with the assumption  $u_{-}>u_{+}$.
Then we must have $\alpha=\beta$, which  means $\lim\limits_{\gamma\rightarrow1}\rho_{\ast}(\gamma)=\alpha.$

If  $\alpha\in(0,+\infty),$ then  we  can also get a contradiction when taking limit in (4.7). Hence $\alpha=0 $ or $ \alpha=+\infty$. By the condition
$\rho_{\ast}>\max\{\rho_{-},\rho_{+}\}$, it is easy to see that $\lim\limits_{\gamma\rightarrow1}\rho_{\ast}(\gamma)=\alpha=+\infty.$

Next  taking the limit $ \gamma\rightarrow 1$ at right hand side in (4.7), we have
$$\lim_{\gamma\rightarrow 1}\sqrt{\frac{\frac{\gamma-1}{2}(\rho_{\ast}^{\gamma-1}-\rho_{\pm}^{\gamma-1})}
{\rho_{\ast}^{2}-\rho_{\pm}^{2}}}(\rho_{\ast}-\rho_{\pm})
=\lim_{\gamma\rightarrow1}\sqrt{\frac{(\frac{\gamma-1}{2}\rho_{\ast}^{\gamma-1}-\frac{\gamma-1}{2}\rho_{\pm}^{\gamma-1})(\rho_{\ast}-\rho_{\pm})^{2}}
{\rho_{\ast}^{2}-\rho_{\pm}^{2}}}
=:\sqrt{a},$$
and
$$u_{-}-u_{+}=2\sqrt{a},\eqno{}
$$
from which we can get
$ a=\frac{(u_{-}-u_{+})^{2}}{4}.$
 The proof is completed. $~\Box$

\vskip 0.1in
\noindent{\small {\small\bf Proposition 4.3.}
It holds that
$$ \lim_{\gamma\rightarrow1}u_{\ast}=\lim_{\gamma\rightarrow1}\sigma_{1}=\lim_{\gamma\rightarrow1}\sigma_{2}=\sigma,\eqno{(4.9)}
$$
and
$$ \lim_{\gamma\rightarrow1}\rho_{\ast}(\sigma_{2}-\sigma_{1})=\sigma[\rho]-[\rho u], \eqno{(4.10)}
$$
where $\sigma=\frac{1}{2}(u_{-}+u_{+}).$

\vskip 0.1in
\noindent{\small {\small\bf Proof.}
Using (3.12), (3.14) and using Lemma 4.2, one  can compute out
$$\lim_{\gamma\rightarrow1}u_{\ast}=u_{-}-\lim_{\gamma\rightarrow1}\sqrt{\frac{\frac{\gamma-1}{2}(\rho_{\ast}^{\gamma-1}-\rho_{-}^{\gamma-1})}
{(\rho_{\ast}+\rho_{-})(\rho_{\ast}-\rho_{-})}}(\rho_{\ast}-\rho_{-})
$$$$
=u_{-}-\sqrt{a}=u_{-}-\frac{1}{2}
(u_{-}-u_{+})=\sigma,
$$
$$\lim_{\gamma\rightarrow1}\sigma_{1}=\lim_{\gamma\rightarrow1}\frac{\rho_{\ast}u_{\ast}-\rho_{-}u_{-}}{\rho_{\ast}-\rho_{-}}
=u_{-}+\lim_{\gamma\rightarrow1}\frac{\rho_{\ast}}{\rho_{-}-\rho_{\ast}}(u_{-}-u_{\ast})
=\sigma,
$$and
$$\lim_{\gamma\rightarrow1}\sigma_{2}=\lim_{\gamma\rightarrow1}\frac{\rho_{+}u_{+}-\rho_{\ast}u_{\ast}}{\rho_{+}-\rho_{\ast}}
=u_{+}+\lim_{\gamma\rightarrow1}\frac{\rho_{\ast}}{\rho_{+}-\rho_{\ast}}(u_{+}-u_{\ast})
=\sigma,
$$which immediately lead to $ \lim\limits_{\gamma\rightarrow1}u_{\ast}=\lim\limits_{\gamma\rightarrow1}\sigma_{1}=\lim\limits_{\gamma\rightarrow1}\sigma_{2}=\sigma.
$

 From the first equations of the Rankine-Hugoniot conditions (3.7) for $S_{1}$ and $S_{2}$,   we have
 $$ \sigma_{1}(\rho_{-}-\rho_{\ast})=\rho_{-}u_{-}-\rho_{\ast}u_{\ast},\eqno{(4.11)}$$and $$\sigma_{2}(\rho_{\ast}-\rho_{+})=\rho_{\ast}u_{\ast}-\rho_{+}u_{+}. \eqno{(4.12)}$$
From  (4.11),
(4.12) and (4.9), we get
 $$ \lim_{\gamma\rightarrow1}\rho_{\ast}(\sigma_{2}-\sigma_{1})=\lim_{\gamma\rightarrow1}(\rho_{-}u_{-}-\sigma_{1}\rho_{-}+\sigma_{2}\rho_{+}-\rho_{+}u_{+})
=\sigma[\rho]-[\rho u].
$$
The proof is completed. $~~\Box$

\vskip 0.1in
\noindent{\small {\small\bf Remark 4.1.} Note that the two shock curves (3.10), (3.11) become very close to the line $u=u_{-}$ as $ \gamma$  tends to one, then
it can be concluded from  Lemma 4.3 that, when $\gamma\rightarrow1$,  the two shock  curves $S_{1}$ and $S_{2}$ will coincide
 to form a new  delta wave,  and   the delta wave   speed $\sigma$ is the limit of both the particle velocity $u_{\ast}$ and
two shocks' speed $\sigma_{1}, \sigma_{2} $.

\vskip 0.1in
What is more, we will further derive that, when $\gamma\rightarrow1$, the limit of Riemann solutions of (1.1) and (1.2) is the delta wave  of the pressureless Euler system (1.4)  with the same Riemann initial data
$(\rho_{\pm}, u_{\pm})$  in the sense of distributions.

\vskip 0.1in
\noindent{\small {\small\bf Theorem 4.4.}
Let $u_{+}<u_{-}.$   For any fixed $\gamma \in(1,2)$, assume that $(\rho_{\gamma}(t,x),u_{\gamma}(t,x))$ is a Riemann solution
containing two shocks $S_{1}$ and $S_{2}$ of (1.1) and (1.2) constructed in Section
3. Then, as $\gamma\rightarrow1$, $(\rho_{\gamma}(t,x),u_{\gamma}(t,x))$ will converge to
 $$(\rho(t,x),u(t,x))=(\rho_{0}(t,x)+w(t)\delta_{S},u_{0}(t,x)),$$
in the sense of distributions, and the singular part of the limit function
$\rho(t, x)$ is  a $\delta$-measure with weight
 $ w(t)=t(\sigma[\rho]-[\rho u])$,  where $\sigma=\frac{1}{2}(u_{-}+u_{+}).$

\vskip 0.1in
\noindent{\small {\small\bf Proof.} (1) Set $\xi=\frac{x}{t}.$  Then for any fixed $\gamma \in(1,2)$, the Riemann solution containing two shocks $S_{1}$ and $S_{2}$ of  (1.1) and  (1.2) can be written as $$(\rho_{\gamma},u_{\gamma})(\xi)=\left\{\begin{array}{ll} (\rho_{-},
u_{-}),\,\,\,\,\xi< \sigma_{1},\\(\rho_{\ast},
u_{\ast}),\,\,\,\,\sigma_{1}< \xi < \sigma_{2},\\(\rho_{+},
u_{+}),\,\,\,\,\xi> \sigma_{2}.\end{array} \right.\eqno{}$$
From (3.3), we have the following weak formulations:
$$-\int_{-\infty}^{+\infty}\rho_{\gamma}(\xi)(u_{\gamma}(\xi)-\xi)\varphi'(\xi)d\xi
+\int_{-\infty}^{+\infty}\rho_{\gamma}(\xi)\varphi(\xi)d\xi=0,
\eqno{(4.13)}$$
$$ \int_{-\infty}^{+\infty}u_{\gamma}(\xi)\varphi(\xi)d\xi
-\int_{-\infty}^{+\infty}\left(\frac{u_{\gamma}(\xi)}{2}-\xi\right)u_{\gamma}(\xi)\varphi'(\xi)d\xi
-\frac{\gamma-1}{4}\int_{-\infty}^{+\infty}\rho_{\gamma}^{\gamma-1}(\xi)\varphi'(\xi)d\xi=0,\eqno{(4.14)}
$$
for any $\varphi(\xi)\in C_{0}^{+\infty}(R)$.

(2) For the second integral on the left-hand side of (4.14),  using the method of integration by parts, we can derive
$$
\int_{-\infty}^{+\infty}(\frac{u_{\gamma}(\xi)}{2}-\xi)u_{\gamma}(\xi)\varphi'(\xi)d\xi
$$
$$=\int_{-\infty}^{\sigma_{1}}(\frac{u_{-}}{2}-\xi)u_{-}\varphi'(\xi)d\xi
+\int_{\sigma_{2}}^{+\infty}(\frac{u_{+}}{2}-\xi)u_{+}\varphi'(\xi)d\xi
+\int_{\sigma_{1}}^{\sigma_{2}}(\frac{u_{\ast}}{2}-\xi)u_{\ast}\varphi'(\xi)d\xi
$$
$$=\frac{u_{-}^{2}}{2}\varphi(\sigma_{1})-\frac{u_{+}^{2}}{2}\varphi(\sigma_{2})
+u_{+}\sigma_{2}\varphi(\sigma_{2})-u_{-}\sigma_{1}\varphi(\sigma_{1})
+u_{-}\int_{-\infty}^{\sigma_{1}}\varphi(\xi)d\xi+u_{+}\int_{\sigma_{2}}^{+\infty}\varphi(\xi)d\xi
+\int_{\sigma_{1}}^{\sigma_{2}}(\frac{u_{\ast}}{2}-\xi)u_{\ast}\varphi'(\xi)d\xi.
$$
Meanwhile, we have
$$\int_{\sigma_{1}}^{\sigma_{2}}(\frac{u_{\ast}}{2}-\xi)u_{\ast}\varphi'(\xi)d\xi
=\frac{u_{\ast}^{2}}{2}(\varphi(\sigma_{2})-\varphi(\sigma_{1}))-u_{\ast}(\sigma_{2}\varphi(\sigma_{2})-\sigma_{1}\varphi(\sigma_{1}))
+u_{\ast}\int_{\sigma_{1}}^{\sigma_{2}}\varphi(\xi)d\xi
$$
$$=u_{\ast}(\sigma_{2}-\sigma_{1})\left( \frac{u_{\ast}}{2}\frac{\varphi(\sigma_{2})-\varphi(\sigma_{1})}{\sigma_{2}-\sigma_{1}}
-\frac{\sigma_{2}\varphi(\sigma_{2})-\sigma_{1}\varphi(\sigma_{1})}{\sigma_{2}-\sigma_{1}}
+\frac{\int_{\sigma_{1}}^{\sigma_{2}}\varphi(\xi)d\xi}{\sigma_{2}-\sigma_{1}}\right).$$
Then, by Proposition 4.3, we have
 $$\lim_{\gamma\rightarrow1}\int_{\sigma_{1}}^{\sigma_{2}}(\frac{u_{\ast}}{2}-\xi)u_{\ast}\varphi'(\xi)d\xi=0.
  $$
Similarly, we can obtain that
$$ \frac{\gamma-1}{4}\int_{-\infty}^{+\infty}\rho_{\gamma}^{\gamma-1}(\xi)\varphi'(\xi)d\xi
=\frac{\gamma-1}{4}\int_{-\infty}^{\sigma_{1}}\rho_{-}^{\gamma-1}\varphi'(\xi)d\xi
+\frac{\gamma-1}{4}\int_{\sigma_{2}}^{+\infty}\rho_{+}^{\gamma-1}\varphi'(\xi)d\xi
+\frac{\gamma-1}{4}\int_{\sigma_{1}}^{\sigma_{2}}\rho_{\ast}^{\gamma-1}\varphi'(\xi)d\xi
 $$
 $$=\frac{\gamma-1}{4}\rho_{-}^{\gamma-1}\varphi(\sigma_{1})-\frac{\gamma-1}{4}\rho_{+}^{\gamma-1}\varphi(\sigma_{2})
 +\frac{\gamma-1}{4}\rho_{\ast}^{\gamma-1}(\varphi(\sigma_{2})-\varphi(\sigma_{1})).
 $$
 It can be derived  from Lemma 4.2 that $ \lim\limits_{\gamma\rightarrow1}\frac{\gamma-1}{2}\rho_{\ast}^{\gamma-1}$  is bounded, then by Proposition 4.3, we have$$\lim_{\gamma\rightarrow1}\frac{\gamma-1}{4}\int_{-\infty}^{+\infty}\rho_{\gamma}^{\gamma-1}(\xi)\varphi'(\xi)d\xi=0.
 $$
 Hence taking the limit $\gamma\rightarrow1$ in (4.14) leads to
$$ \lim\limits_{\gamma\rightarrow1}\int_{-\infty}^{+\infty}(u_{\gamma}(\xi)-u_{0}(\xi))\varphi(\xi)d\xi
=\left(\sigma[u]-[\frac{u^{2}}{2}]\right)\varphi(\sigma)=\left(\frac{1}{2}(u_{-}+u_{+})[u]-[\frac{u^{2}}{2}]\right)\varphi(\sigma)=0,\eqno (4.15)
$$
where $ (\rho_{0}(\xi),u_{0}(\xi))=(\rho_{\pm},u_{\pm}),~\pm(\xi-\sigma)>0.$

(3) Similarly, we can obtain for (4.13) that
$$\lim_{\gamma\rightarrow1}\int_{-\infty}^{+\infty}(\rho_{\gamma}(\xi)-\rho_{0}(\xi))\varphi(\xi)d\xi
=(\sigma[\rho]-[\rho u])\varphi(\sigma).\eqno (4.16)
$$

(4) Finally, we study the limits of $\rho_{\gamma}(t,x) $ and $u_{\gamma}(t,x) $ depending on $t$ as $\gamma\rightarrow1$. For any  $ \varphi(t,x)\in C_{0}^{+\infty}(R_{+}\times R)$,    we have
$$\lim_{\gamma\rightarrow1}\int_{0}^{+\infty}\int_{-\infty}^{+\infty}\rho_{\gamma}(x/t)\varphi(t,x)dx dt
=\lim_{\gamma\rightarrow1}\int_{0}^{+\infty}t\bigg(\int_{-\infty}^{+\infty}\rho_{\gamma}(\xi)\varphi(t,\xi t)d\xi\bigg) dt.
\eqno (4.17)
$$ Regarding $t$ as a parameter, one can  get from   (4.16) that
$$
\lim_{\gamma\rightarrow1}\int_{-\infty}^{+\infty}\rho_{\gamma}(\xi)\varphi(t,\xi t)d\xi=\int_{-\infty}^{+\infty}\rho_{0}(\xi)\varphi(t,\xi t)d\xi+(\sigma[\rho]-[\rho u])\varphi(t,\sigma t)$$$$=\frac{1}{t}\int_{-\infty}^{+\infty}\rho_{0}(x/t)\varphi(t,x)dx+(\sigma[\rho]-[\rho u])\varphi(t,\sigma t)$$$$=\frac{1}{t}\int_{-\infty}^{+\infty}\rho_{0}(t,x)\varphi(t,x)dx+(\sigma[\rho]-[\rho u])\varphi(t, \sigma t).\eqno (4.18)
$$ Substituting (4.18) into (4.17),  we have
$$\lim_{\gamma\rightarrow1}\int_{0}^{+\infty}\int_{-\infty}^{+\infty}\rho_{\gamma}(x/t)\varphi(t,x)dxdt
=\int_{0}^{+\infty}\int_{-\infty}^{+\infty}\rho_{0}(t,x)\varphi(t,x)dxdt+\int_{0}^{+\infty}t(\sigma[\rho]-[\rho u])\varphi(t,\sigma t)dt.
$$
This yields that
$$\lim_{\gamma\rightarrow1}\int_{0}^{+\infty}\int_{-\infty}^{+\infty}(\rho_{\gamma}(t,x)-\rho_{0}(t,x))\varphi(t,x)dx dt
$$ $$
=\int_{0}^{+\infty}t(\sigma[\rho]-[\rho u])\varphi(t,\sigma t)dt,\eqno {(4.19)}
$$
in which by definition (2.1), we have
$$\int_{0}^{+\infty}t(\sigma[\rho]-[\rho u])\varphi(t,\sigma t)dt=
\langle w(\cdot)\delta_{S},\varphi(\cdot,\cdot)\rangle.\eqno (4.20)
$$
where
 $$ w(t)=t(\sigma[\rho]-[\rho u]).$$

Similar to (4.18), we can derive from (4.15) that
$$
\lim_{\gamma\rightarrow1}\int_{-\infty}^{+\infty}u_{\gamma}(\xi)\varphi(t,\xi t)d\xi=\int_{-\infty}^{+\infty}u_{0}(\xi)\varphi(t,\xi t)d\xi=\frac{1}{t}\int_{-\infty}^{+\infty}u_{0}(x/t)\varphi(t,x)dx=\frac{1}{t}\int_{-\infty}^{+\infty}u_{0}(t,x)\varphi(t,x)dx.\eqno (4.21)
$$
Therefore,$$\lim_{\gamma\rightarrow1}\int_{0}^{+\infty}\int_{-\infty}^{+\infty}u_{\gamma}(x/t)\varphi(t,x)dx dt
=\lim_{\gamma\rightarrow1}\int_{0}^{+\infty}t\bigg(\int_{-\infty}^{+\infty}u_{\gamma}(\xi)\varphi(t,\xi t)d\xi\bigg) dt$$
$$
=\int_{0}^{+\infty}\int_{-\infty}^{+\infty}u_{0}(t,x)\varphi(t,x)dxdt,
$$
which implies that$$\lim_{\gamma\rightarrow1}\int_{0}^{+\infty}\int_{-\infty}^{+\infty}(u_{\gamma}(t, x)-u_{0}(t,x))\varphi(t,x)dxdt
=0.
$$
The proof is completed. $~~\Box$

\baselineskip 15pt
 \sec{\Large\bf 5.\quad  Numerical results}

 In order to verify the validity of the formation  the delta wave mentioned in section
  4, in this section we present  a selected group of representative numerical results
by using Euler system (1.1), with the Riemann initial data (1.2). A number of iterative numerical trials are executed to guarantee what we demonstrate are not numerical objects. To discretize the system (1.1), we use the
  fifth-order weighted essentially non-oscillatory scheme and third-order Runge-Kutta method  [24, 25]  with the mesh 200 points.

When $u_{+}<u_{-}$, we compute the solution of the Riemann problem of (1.1)-(1.2)  and take
 the initial data  as follows:
$$ (\rho, u)(0, x) =\left\{\begin{array}{ll} (1.5,
1.5),\,\,\,\,x< 0,\\(2,
-0.5),\,\,\,\,x> 0.\end{array} \right.\eqno{(5.1)}$$
The numerical simulations for different choices of $\gamma$ (i.e., $\gamma$=2.5,~1.3,~1.05,~1.0001, and the time $t=0.3$)  are presented in Figs. 2-5 which show the process of concentration of mass and formation of the delta wave  in the pressureless limit of solutions containing two shocks.

\begin{figure}[htbp]
\centering
\begin{minipage}[c]{0.45\textwidth}
\centering
\includegraphics[width=2.5in]{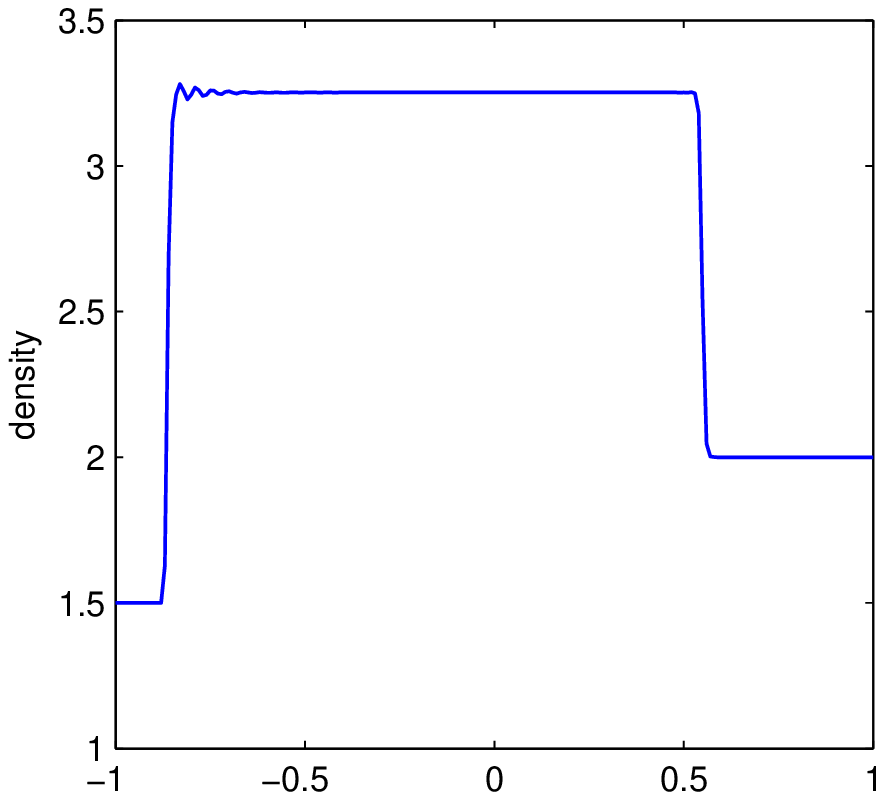}
\end{minipage}%
\begin{minipage}[c]{0.45\textwidth}
\centering
\includegraphics[width=2.5in]{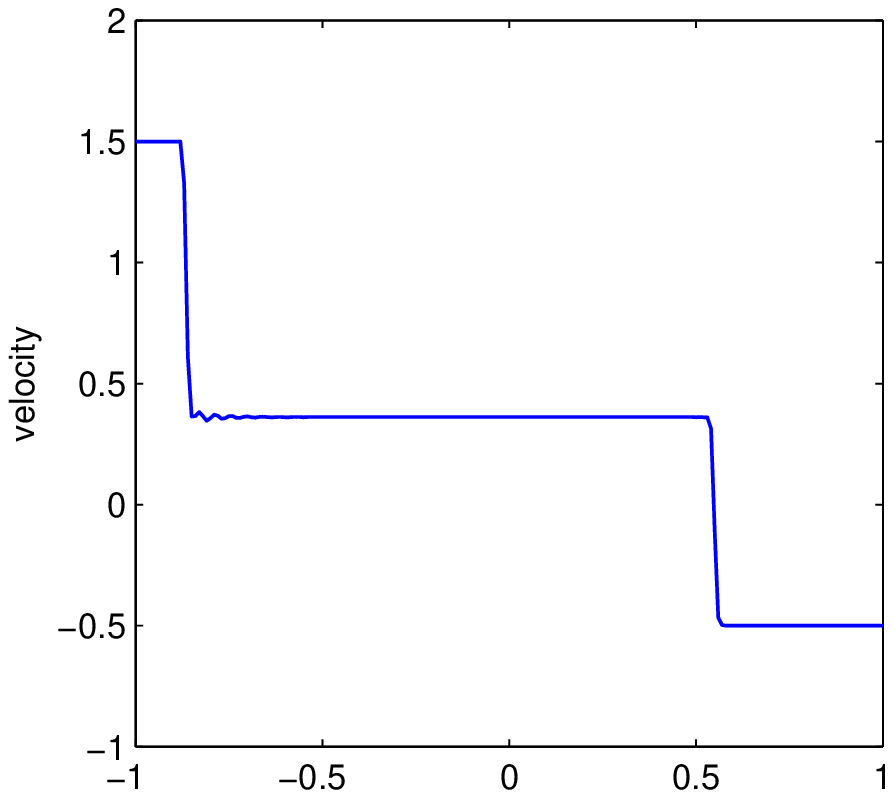}
\end{minipage}%
\end{figure}

\centerline{\bf Fig. 2.\,\,   Density (left) and velocity (right) for $\gamma=2.5$.}

\begin{figure}[htbp]
\centering
\begin{minipage}[c]{0.45\textwidth}
\centering
\includegraphics[width=2.5in]{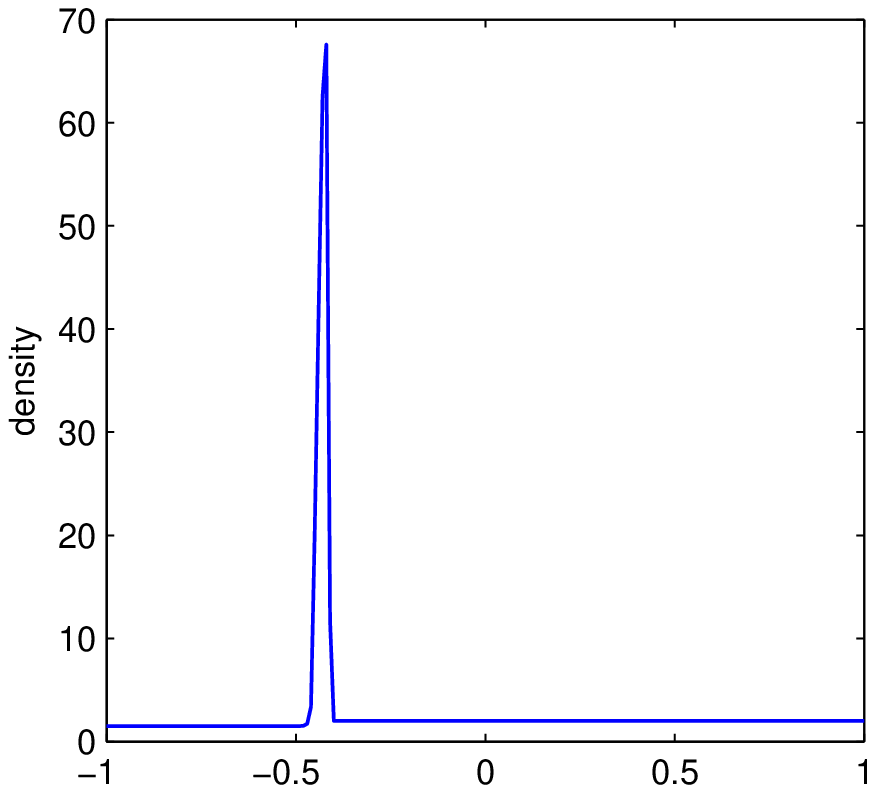}
\end{minipage}%
\begin{minipage}[c]{0.45\textwidth}
\centering
\includegraphics[width=2.5in]{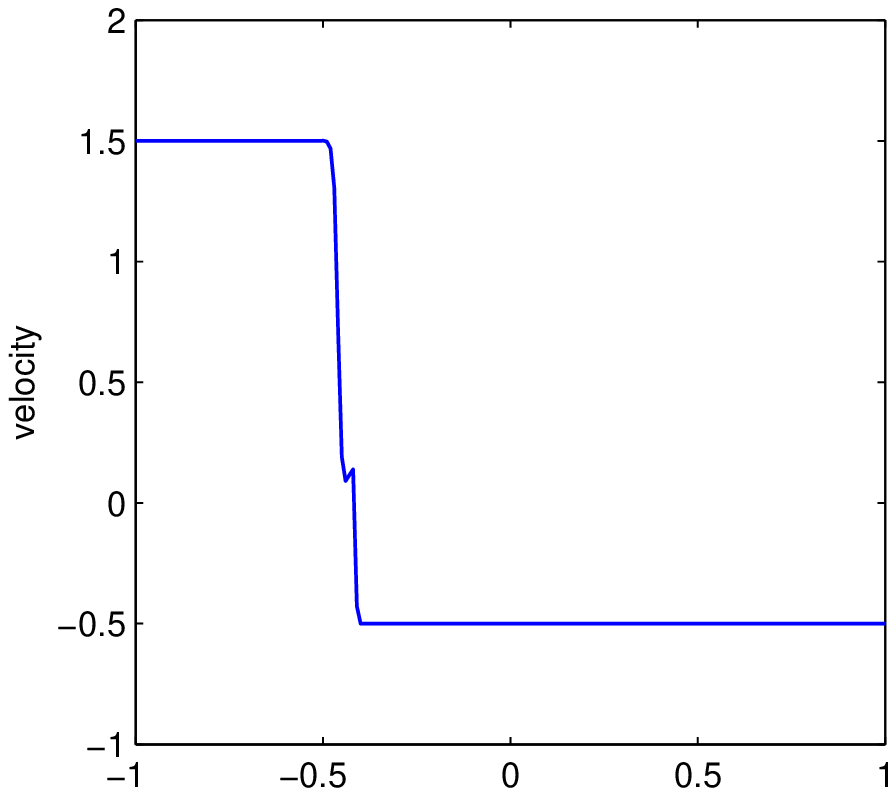}
\end{minipage}%
\end{figure}

\centerline{\bf Fig. 3.\,\,   Density (left) and velocity (right) for $\gamma=1.3$.}

\newpage
\begin{figure}[htbp]
\centering
\begin{minipage}[c]{0.45\textwidth}
\centering
\includegraphics[width=2.5in]{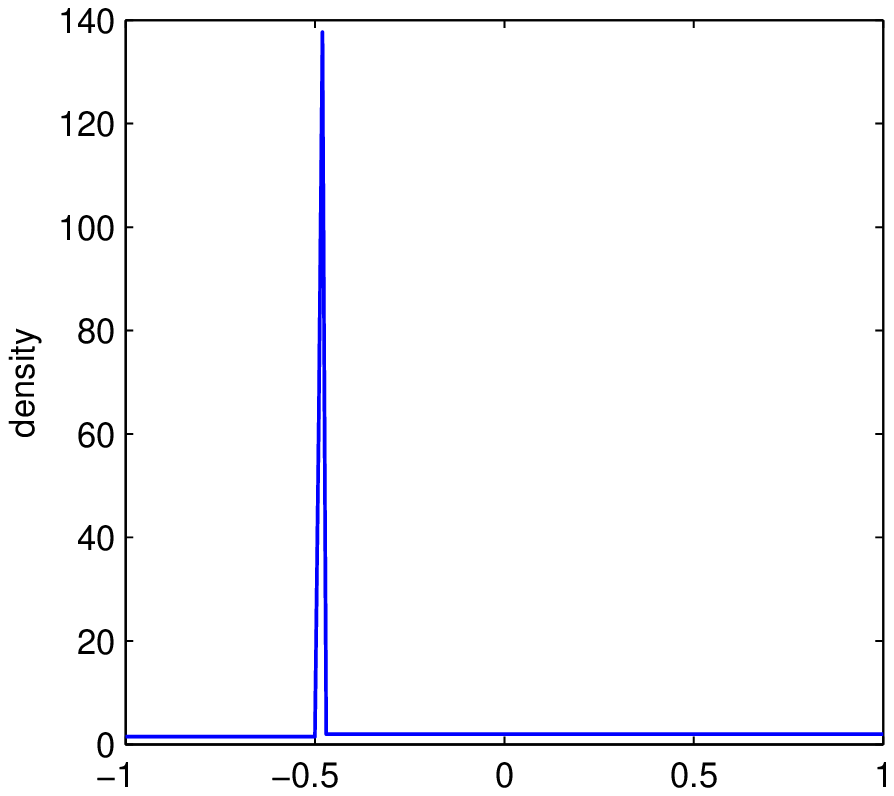}
\end{minipage}%
\begin{minipage}[c]{0.45\textwidth}
\centering
\includegraphics[width=2.5in]{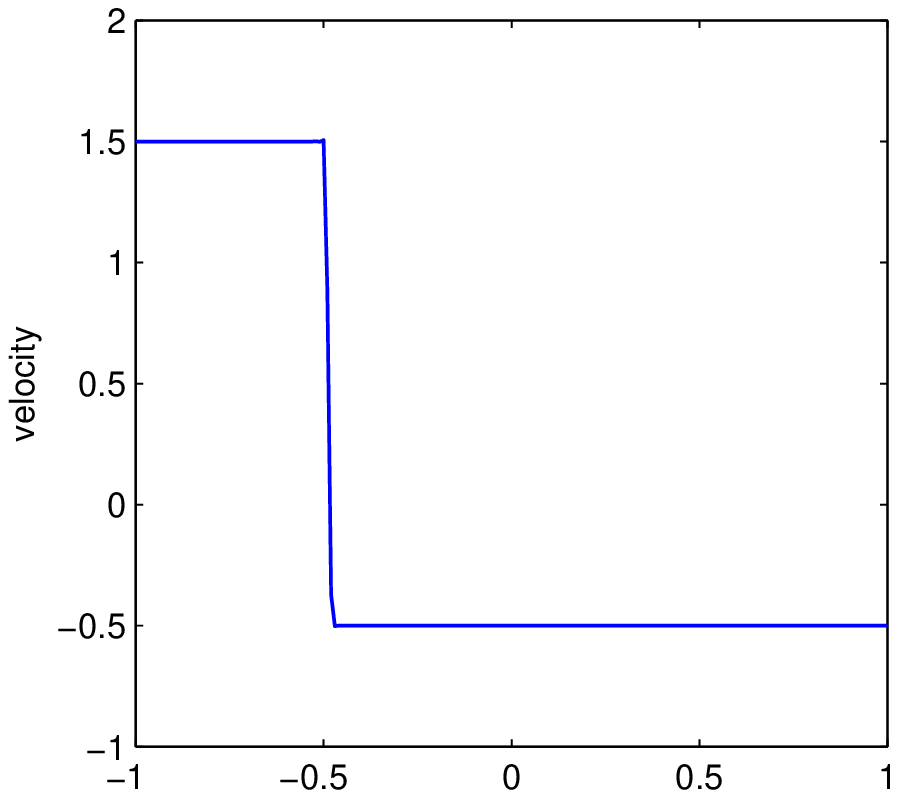}
\end{minipage}%
\end{figure}

\centerline{\bf Fig. 4.\,\,   Density (left) and velocity (right) for $\gamma=1.05$.}

\begin{figure}[htbp]
\centering
\begin{minipage}[c]{0.45\textwidth}
\centering
\includegraphics[width=2.5in]{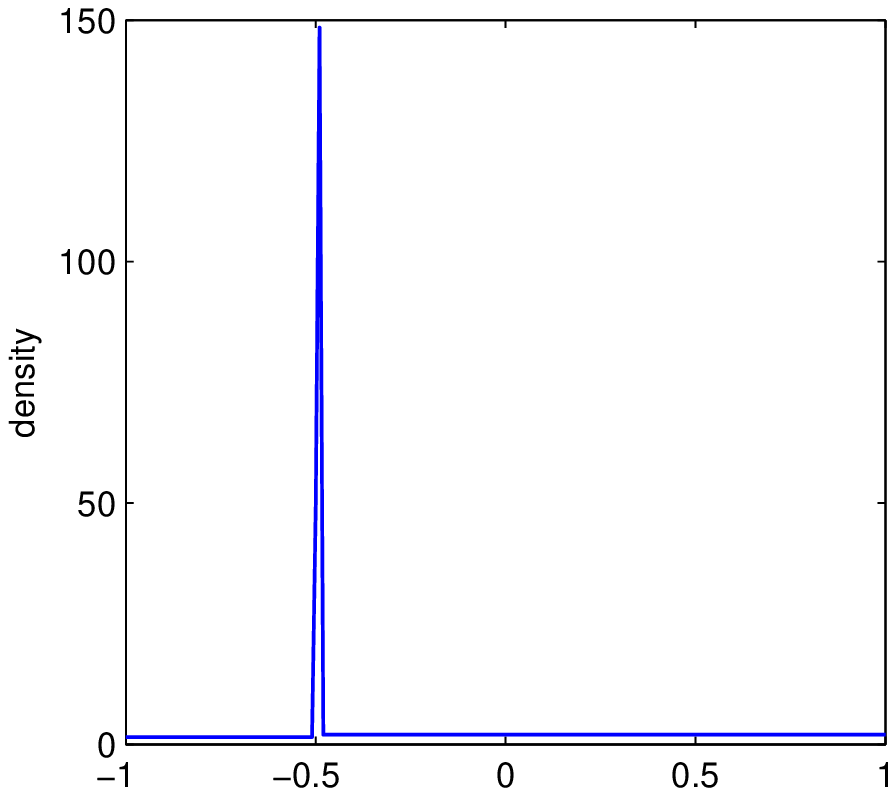}
\end{minipage}%
\begin{minipage}[c]{0.45\textwidth}
\centering
\includegraphics[width=2.5in]{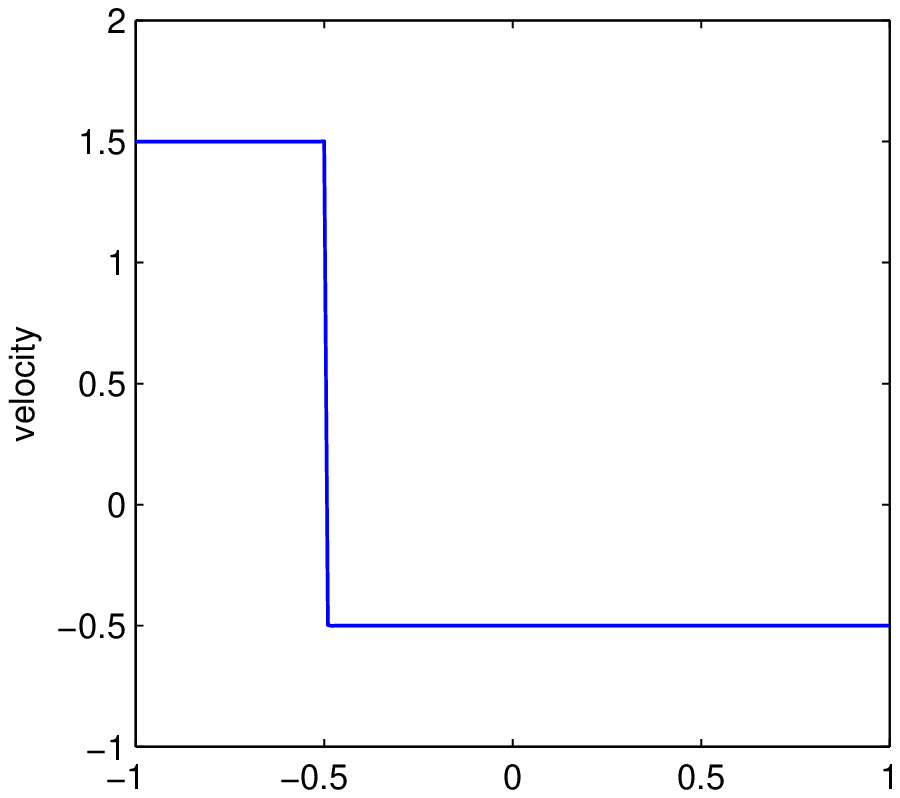}
\end{minipage}%
\end{figure}

\centerline{\bf Fig. 5.\,\,   Density (left) and velocity (right) for $\gamma=1.0001$.}

We can clearly observe  from these numerical results  that,  as $\gamma$ decreases, the locations of the two shocks become closer, and the density of the intermediate state increases dramatically, while  the velocity becomes a piecewise constant function. In the end, as $\gamma\rightarrow1$, along with the intermediate state, the two shocks coincide to form  the delta wave of the pressureless
Euler system (1.4), (1.2), while the velocity is a piecewise constant function. The numerical simulations are in complete agreement with
 the theoretical analysis in section 4.

 \indent

\end{document}